\documentclass[11pt]{article}
\usepackage{amsfonts}
\usepackage{amsfonts}
\usepackage{amsfonts}
\usepackage{mathrsfs}
\usepackage{diagbox}
\usepackage{bm}
\usepackage{cite}
\usepackage[colorlinks,linkcolor=blue,citecolor=blue]{hyperref}
\usepackage{amssymb,amsmath} 
\usepackage{booktabs}
 \allowdisplaybreaks

\catcode`!=11
\let\!int\int \def\int{\displaystyle\!int}
\let\!lim\lim \def\lim{\displaystyle\!lim}
\let\!sum\sum \def\sum{\displaystyle\!sum}
\let\!sup\sup \def\sup{\displaystyle\!sup}
\let\!inf\inf \def\inf{\displaystyle\!inf}
\let\!cap\cap \def\cap{\displaystyle\!cap}
\let\!max\max \def\max{\displaystyle\!max}
\let\!min\min \def\min{\displaystyle\!min}
\let\!frac\frac \def\frac{\displaystyle\!frac}
\catcode`!=12

\let\oldsection\section
\renewcommand\section{\setcounter{equation}{0}\oldsection}

\allowdisplaybreaks
\def\pf{\it{Proof.}\rm\quad}

\def\R{\mathbb{R}}
\def\N{\mathbb{N}}\def\Z{\mathbb{Z}}

\newtheorem{thm}{Theorem}[section]
\newtheorem{lem}[thm]{Lemma}
\newtheorem{cor}[thm]{Corollary}

\def\su{\sum\limits_{n=1}^\infty}
\begin{document}
\title {\bf Identities about level 2 Eisenstein series}
\author{
{Ce Xu\thanks{Corresponding author. Email: xuce1242063253@163.com}}\\[1mm]
\small School of Mathematical Sciences, Xiamen University\\
\small Xiamen
361005, P.R. China}

\date{}
\maketitle \noindent{\bf Abstract } In this paper we consider certain classes of generalized level 2 Eisenstein series by simple differential calculations of trigonometric functions. In particular, we give four new transformation formulas for some level 2 Eisenstein series. We can find that these level 2 Eisenstein series are reducible to infinite series involving hyperbolic functions. Moreover, some interesting new examples are given.
\\[2mm]
\\[2mm]

\section{Introduction}

Let $\N$ be the set of natural numbers, $\N_0:=\N\cup \{0\}$, $\Z$ the ring of integers, $\mathbb{Q}$ the field of rational numbers, $\R$ the field of real numbers, and $\mathbb{C}$ the field of complex numbers. Let $i=\sqrt{-1}$.

The subject of this paper are Eisenstein series and hyperbolic functions. Let $\tau$ be a complex number with strictly positive imaginary part, the holomorphic Eisenstein series $G_{2k}(\tau)$ of weight $2k$, where $k\geq2$ is an integer, is defined by the following series:
\begin{align}
G_{2k}(\tau):=\sum\limits_{m,n\in \Z\setminus(0,0)}\frac{1}{(m+n\tau)^{2k}}.
\end{align}
This series absolutely converges to a holomorphic function of $\tau$ in the upper half-plane. It is well known that the value of $G_{4k}(i)$ can be expressed as
\begin{align}
G_{4k}(i)=\frac{\Gamma^{8k}(1/4)}{2^{2k}\pi^{4k}}H_{4k}\quad (k\in\N),
\end{align}
where $H_{4m}$ are called the Hurwitz numbers (see \cite{A1974,H1962,T2008}).
When working with the $q$-expansion of the Eisenstein series, this alternate notation is frequently introduced:
\begin{align}\label{1.2}
E_{2k}(\tau):=\frac{G_{2k}(\tau)}{2\zeta(2k)}=1+\frac{2}{\zeta(1-2k)}\sum\limits_{n=1}^\infty\frac{n^{2k-1}q^n}{1-q^n},
\end{align}
where $q=\exp(2\pi \tau)$, and $\zeta(s)$ denotes the Riemann zeta function. In Ramanujan's notation, the three relevant Eisenstein series are defined for $|q|<1$ by
\begin{align}
&P(q):=1-24\su \frac{nq^n}{1-q^n},\\
&Q(q):=1+240\su \frac{n^3q^n}{1-q^n},\\
&R(q):=1-504\su \frac{n^5q^n}{1-q^n}.
\end{align}
Thus, for $q=\exp(2\pi i \tau)$, $E_4(\tau)=Q(q)$ and $E_6(\tau)=R(q)$, which have weights 4 and 6, respectively. Since (\ref{1.2}) does not coverge for $j=1$, the Eisenstein series $E_2(\tau)$ must be defined differently, which is defied by
\begin{align}
E_2(\tau)=P(q)-\frac{3}{\pi {\rm Im} \tau}.
\end{align}
The functions $P, Q$ and $R$ were thoroughly studied in a famous paper \cite{H1927} by Ramanujan. Berndt \cite{B1977,B1978} found a lot identities about infinite series involving hyperbolic functions using certain modular transformation formula that originally stems from the gereralized Eisensein series. Further results of infinite series involving hyperbolic functions see Berndt's books \cite{B1989,B1991} and the references therein.

Recently, surprisingly little work has been done on level 2 Eisenstein series involving hyperbolic functions. The motivation for this paper arises from results of Tsumura \cite{T2008,T2009,T2010,T2012,T2018} and \cite{T2015} with Komori and Matsumoto. They studied many level 2 Eisenstein series involving hyperbolic functions. For example, in 2008, Tsumura \cite{T2008} considered the following two level 2 Eisenstein series of hyperbolic functions
\begin{align*}
\mathcal{G}_k(i):=\sum\limits_{m,n\in\Z,m\neq 0} \frac{(-1)^n}{\sinh(m\pi)(m+ni)^k}\quad {\rm and}\quad \mathcal{H}_k(i):=\sum\limits_{m,n\in\Z,m\neq 0} \frac{(-1)^n}{\cosh(m\pi)(m+ni)^k}\quad (k\in\N).
\end{align*}
He proved that $\mathcal{G}_{2k-1}(i)$ and $\mathcal{H}_{2k}(i)$ can be expressed in terms of $\Gamma$ function and $\pi$. Further, in 2009, Tsumura \cite{T2009} studied the closed form representations of sums
\begin{align*}
\mathcal{C}_k^v:=\sum\limits_{m,n\in\Z,m\neq 0} \frac{\coth^v(m\pi)}{(m+ni)^k}
\end{align*}
for $k\in\N$ with $k\geq 3$ and $v\in\Z$. Specially, he showed that
\begin{align*}
\mathcal{C}_k^v\in \mathbb{Q}\left[\frac 1{\pi},\pi,\frac{\Gamma^8(1/4)}{\pi^2}\right]
\end{align*}
for $k\geq 3$ and $v\in\N_0$ with $k\equiv v$ (mod 2).

In this paper, continuing Tsumura et al's work, we study the four level 2 Eisenstein series
\begin{align*}
\sum\limits_{m,n\in\Z,\atop m\neq 0} \frac{f(m)(-1)^n}{(m+ani)^k},\ \sum\limits_{m,n\in\Z,\atop m\neq 0} \frac{f(m)(-1)^n}{(m+a(2n+1)i)^k},\
\sum\limits_{m,n\in\Z,\atop m\neq 0} \frac{g(m)}{(m+ani)^p},\ \sum\limits_{m,n\in\Z,\atop m\neq 0} \frac{g(m)}{(m+a(2n+1)i)^p},
\end{align*}
 where $k\in\N,\ p\in\N\setminus \{1\},\ a\in \mathbb{R}\setminus \{0\}$ and if $\quad m\rightarrow \infty$
\begin{align*}
&f(m)=o(1),\\
&g(m)=o(1/m),\quad p=2 \quad{\rm and}\quad g(m)=o(1), \quad p>2.
\end{align*}
We prove that these double series can be expressed by single infinite series involving hyperbolic functions. Moreover, we consider some special cases. We can find that many level 2 Eisenstein series involving hyperbolic functions can be expressed in terms of $\Gamma$ function and $\pi$.

\section{Differential formulas of trigonometric functions}

Let \[|{\bf r}|_{l}:=r_0+r_1+\cdots+r_{l} \ (r_j\in\N_0).\]

\begin{lem} Let $I_{k,m}$ be a sequence and $k$ and $m$ are positive integers (including zero). If $I_{k,m}$ satisfy a recurrence relation in the form
\begin{align}\label{2.1}
I_{m,k}=a_mI_{k-1,m+1}-b_mI_{k-1,m}
\end{align}
then
\begin{align}\label{2.2}
I_{k,m}=\sum\limits_{l=0}^{k-1}\left(\prod\limits_{j=0}^l a_{m+j}\right)\left(\sum\limits_{|{\bf r}|_{l+1}=k-l-1} \prod\limits_{h=0}^{l+1}b_{m+h}^{r_{h}}\right)I_{0,m+l+1}(-1)^{k-l+1}+(-1)^kb_m^kI_{0,m},
\end{align}
where $a_m$ and $b_m$ are constants.
\end{lem}
\pf The result (\ref{2.2}) can be proved by mathematical induction. \hfill$\square$

\begin{thm}\label{thm2.2} For integers $k\geq 0,\ m\geq 1$ and complex number $s\in \mathbb{C}\setminus \N_0$, then
\begin{align}
&\frac{d^{2k}}{ds^{2k}}\left(\frac{1}{\sin^{2m-1}(\pi s)}\right)=\pi^{2k} \sum\limits_{l=0}^k \left(\frac{(2m+2l-2)!}{(2m-2)!}\sum\limits_{\mid {\bf r} \mid_l=k-l}\prod\limits_{h=0}^l (2m+2h-1)^{2r_h}\right)\frac{(-1)^{k-l}}{\sin^{2m+2l-1}(\pi s)},\label{2.3}\\
&\frac{d^{2k+1}}{ds^{2k+1}}\left(\frac{1}{\sin^{2m-1}(\pi s)}\right)=\pi^{2k+1} \sum\limits_{l=0}^k \left(\frac{(2m+2l-1)!}{(2m-2)!}\sum\limits_{\mid {\bf r} \mid_l=k-l}\prod\limits_{h=0}^l (2m+2h-1)^{2r_h}\right)\nonumber\\&\quad\quad\quad\quad\quad\quad\quad\quad\quad\quad\quad\quad\quad\quad\quad\quad\times\frac{(-1)^{k+1-l}\cos(\pi s)}{\sin^{2m+2l}(\pi s)}.\label{2.4}
\end{align}
\end{thm}
\pf A elementary calculation gives
\begin{align}
\frac{d^{2k}}{ds^{2k}}\left(\frac{1}{\sin^{2m-1}(\pi s)}\right)=&2m(2m-1)\pi^2 \frac{d^{2k-2}}{ds^{2k-2}}\left(\frac{1}{\sin^{2m+1}(\pi s)}\right)\nonumber\\&-(2m-1)^2\pi^2 \frac{d^{2k-2}}{ds^{2k-2}}\left(\frac{1}{\sin^{2m-1}(\pi s)}\right).\label{2.5}
\end{align}
So, setting $I_{k,m}=d^{2k}(1/\sin^{2m-1}(\pi s))/ds^{2k},\ a_m=2m(2m-1)\pi^2$ and $b_m=(2m-1)^2\pi^2$ in (\ref{2.2}) and combining (\ref{2.5}) yields the desired result (\ref{2.3}). Then, differentiating (\ref{2.3}) with respect to $s$, we may deduce the evaluation (\ref{2.4}).\hfill$\square$

In (\ref{2.3}), let
\begin{align}\label{e1}
A_{k,m}(l):=(-1)^{k-l}\frac{(2m+2l-2)!}{(2m-2)!}\sum\limits_{\mid {\bf r} \mid_l=k-l}\prod\limits_{h=0}^l (2m+2h-1)^{2r_h},
\end{align}
we have
\begin{align}\label{e2}
\frac{d^{2k}}{ds^{2k}}\left(\frac{1}{\sin^{2m-1}(\pi s)}\right)=\pi^{2k} \sum\limits_{l=0}^k \frac{A_{k,m}(l)}{\sin^{2m+2l-1}(\pi s)}.
\end{align}
To evaluate $A_{k,m}(l)$, we differentiate both sides of either equation in (\ref{e2}) twice and equate the coefficients. The following recurrence relation is
obtained, for $1\leq l \leq k-1$ and $m\geq 1$:
\begin{align}
A_{k,m}(l)=(2m+2l-3)(2m+2l-2) A_{k-1,m}(l-1)-(2m+2l-1)^2A_{k-1,m}(l).
\end{align}
In particular, for $k\geq 0$ and $m\geq 1$,
\begin{align*}
A_{k,m}(0)=(-1)^k (2m-1)^{2k}\quad {\rm and}\quad A_{k,m}(k)=\frac{(2m+2k-2)!}{(2m-2)!}.
\end{align*}
Hence, in below, we have
\begin{align*}
\sum\limits_{\mid{\bf r}\mid_l=k-l}\prod\limits_{h=0}^l (2h+1)^{2r_h}=(-1)^{k-l}\frac{A_{k,1}(l)}{(2l)!}.
\end{align*}

We give some values of coefficient $A_{k,1}(l)$ with the help of Mathematica.
\begin{table}[!htbp]
\centering
\begin{tabular}{|c|c|c|c|c|c|c|c|c|}
\hline
\diagbox{$l$}{$A_{k,1}(l)$}{$k$}&$0$&$1$&$2$&$3$ & $4$& $5$& $6$& $7$\\ 
\hline
$0$&1 &-1 &1 & -1&1¡¡ &-1¡¡&1¡¡&-1¡¡\\
\hline
$1$&0 &2&-20&182&¡¡-1640 &14762¡¡&-132860¡¡&1195742¡¡\\
\hline
$2$&0 &0&24&-840&23184 &-599280¡¡&¡¡15159144&-380572920¡¡\\
\hline
$3$&0 &0&0&720&-60480 &3659040¡¡&-197271360¡¡&10121070960¡¡\\
\hline
$4$&0 &0&0&0& 40320&¡¡-6652800&743783040¡¡&-71293622400¡¡\\
\hline
$5$&0 &0&0&0& 0&362880¡¡&-1037836800¡¡&192518726400¡¡\\
\hline
\end{tabular}
\\
 [2mm]
  \centering
{\bf {\small Table 1}. coefficient $A_{k,1}(l)$}
\end{table}

\begin{thm}\label{thm2.3} For integer $k\in\N_0$ and complex number $s\in \mathbb{C}\setminus \N_0$, we have
\begin{align}
\frac{d^{2k}}{ds^{2k}}\left(\cot(\pi s)\right)=&\pi^{2k}\sum\limits_{0\leq l\leq j \leq k} \left\{ \binom{2k}{2j}(2l)!-\binom{2k}{2j+1}(2l+1)! \right\}\nonumber\\&\quad\quad\quad\quad\times \left\{\sum\limits_{\mid{\bf r}\mid_l=j-l}\prod\limits_{h=0}^l (2h+1)^{2r_h}\right\}\frac{(-1)^{k-l}\cos(\pi s)}{\sin^{2l+1}(\pi s)},\label{2.6}\\
\frac{d^{2k+1}}{ds^{2k+1}}\left(\cot(\pi s)\right)=&\pi^{2k+1}\sum\limits_{1\leq l\leq j \leq k+1} \left\{ \binom{2k+2}{2j}(2l-1)!-\binom{2k+2}{2j+1}(2l-1)!(2l+1) \right\}\nonumber\\&\quad\quad\quad\quad\quad\times \left\{\sum\limits_{\mid{\bf r}\mid_l=j-l}\prod\limits_{h=0}^l (2h+1)^{2r_h}\right\}\frac{(-1)^{k-l}}{\sin^{2l}(\pi s)},\label{2.7}
\end{align}
where $$\binom{n}{k}:=\frac{n!}{k!(n-k)!},$$
if $k>n$, then $\binom{n}{k}=0$.
\end{thm}
\pf By a direct calculation we find that
\begin{align}
\frac{d^{2k}}{ds^{2k}}\left(\cot(\pi s)\right)&=\frac{d^{2k}}{ds^{2k}}\left(\frac{\sin(\pi s)}{\cos(\pi s)}\right)\nonumber\\
&=\sum\limits_{j=0}^k \binom{2k}{2j}\frac{d^{2j}}{ds^{2j}}\left(\frac{1}{\sin(\pi s)}\right)\frac{d^{2k-2j}}{ds^{2k-2j}}\left(\cos(\pi s)\right)\nonumber\\
&\quad-\sum\limits_{j=0}^k \binom{2k}{2j+1}\frac{d^{2j+1}}{ds^{2j+1}}\left(\frac{1}{\sin(\pi s)}\right)\frac{d^{2k-2j-1}}{ds^{2k-2j-1}}\left(\cos(\pi s)\right).\label{2.8}
\end{align}
Hence, letting $k=j$ and $m=1$ in (\ref{2.3}) and (\ref{2.4}), then substituting it into (\ref{2.8}) we obtain (\ref{2.6}). Integrating (\ref{2.6}) over the interval $(1/2,s)$ with respect to $s$, a simple calculation gives the formula (\ref{2.7}).\hfill$\square$

Further, changing $s$ to $1/2-s$ in Theorems \ref{thm2.2} and \ref{thm2.3}, we can get the following corollaries.
\begin{cor}\label{cor2.4} For integer $k$ and complex number $s$ with $s\neq \pm 1/2, \pm 3/2, \cdots$, then
\begin{align}
&\frac{d^{2k}}{ds^{2k}}\left(\frac {1}{\cos(\pi s)}\right)=\pi^{2k} \sum\limits_{l=0}^k \left((2l)!\sum\limits_{\mid {\bf r} \mid_l=k-l}\prod\limits_{h=0}^l (2h+1)^{2r_h}\right)\frac{(-1)^{k-l}}{\cos^{2l+1}(\pi s)},\label{2.9}\\
&\frac{d^{2k+1}}{ds^{2k+1}}\left(\frac {1}{\cos(\pi s)}\right)=\pi^{2k+1} \sum\limits_{l=0}^k \left((2l+1)!\sum\limits_{\mid {\bf r} \mid_l=k-l}\prod\limits_{h=0}^l (2h+1)^{2r_h}\right)\frac{(-1)^{k-l}\sin(\pi s)}{\cos^{2l+2}(\pi s)}.\label{2.10}
\end{align}
\end{cor}

\begin{cor}\label{cor2.5} For integer $k$ and complex number $s$ with $s\neq \pm 1/2, \pm 3/2, \cdots$, then
\begin{align}
\frac{d^{2k}}{ds^{2k}}\left(\tan(\pi s)\right)=&\pi^{2k}\sum\limits_{0\leq l\leq j \leq k} \left\{ \binom{2k}{2j}(2l)!-\binom{2k}{2j+1}(2l+1)! \right\}\nonumber\\&\quad\quad\quad\quad\times\left\{ \sum\limits_{\mid{\bf r}\mid_l=j-l}\prod\limits_{h=0}^l (2h+1)^{2r_h}\right\}\frac{(-1)^{k-l}\sin(\pi s)}{\cos^{2l+1}(\pi s)},\label{2.11}\\
\frac{d^{2k+1}}{ds^{2k+1}}\left(\tan(\pi s)\right)=&\pi^{2k+1}\sum\limits_{1\leq l\leq j \leq k+1} \left\{ \binom{2k+2}{2j}(2l-1)!-\binom{2k+2}{2j+1}(2l-1)!(2l+1) \right\}\nonumber\\&\quad\quad\quad\quad\quad\times \left\{\sum\limits_{\mid{\bf r}\mid_l=j-l}\prod\limits_{h=0}^l (2h+1)^{2r_h}\right\}\frac{(-1)^{k-l+1}}{\cos^{2l}(\pi s)}.\label{2.12}
\end{align}
\end{cor}

\section{Main Theorems and Corollaries}

In this section we consider the following level 2 Eisenstein series

\begin{align*}
\sum\limits_{m,n\in\Z,m\neq 0} \frac{f(m)}{(m+ani)^k}(-1)^n,\quad \sum\limits_{m,n\in\Z,m\neq 0} \frac{f(m)}{(m+a(2n+1)i)^k}(-1)^n
\end{align*}
and
\begin{align*}
\sum\limits_{m,n\in\Z,m\neq 0} \frac{g(m)}{(m+ani)^p},\quad \sum\limits_{m,n\in\Z,m\neq 0} \frac{g(m)}{(m+a(2n+1)i)^p},
\end{align*}
 where $k\in\N,\ p\in\N\setminus \{1\},\ a\in \mathbb{R}\setminus \{0\}$ and
\begin{align*}
f(m)=o(1),\ g(m)=o(1/m),\quad m\rightarrow \infty.
\end{align*}
Note that if $k=1$ in the first two sums, then
\begin{align*}
&\sum\limits_{m,n\in\Z,m\neq 0} \frac{f(m)(-1)^n}{m+ani}=
\sum\limits_{m\in\Z,m\neq 0}\lim_{N\to\infty}\sum\limits_{-N\leq n\leq N} \frac{f(m)(-1)^n}{m+ani},\\
&\sum\limits_{m,n\in\Z,m\neq 0} \frac{f(m)(-1)^n}{m+a(2n+1)i}=\sum\limits_{m\in\Z,m\neq 0} \lim_{N\to\infty}\sum\limits_{-N\leq n\leq N} \frac{f(m)(-1)^n}{m+a(2n+1)i}.
\end{align*}

\subsection{Four Theorems}

According to the partial fraction expansion of trigonometric function
\begin{align*}
&\frac{\pi}{\sin(\pi s)}=\sum\limits_{n\in\Z}\frac{(-1)^n}{n+s},\quad \frac{\pi}{\cos(\pi s)}=2\sum\limits_{n\in\Z} \frac{(-1)^n}{2n+1-2s},\\
&\pi\cot(\pi s)=\lim_{N\rightarrow\infty}\sum\limits_{-N\leq n\leq N}\frac{1}{n+s},\quad \pi\tan(\pi s)=2\lim_{N\rightarrow\infty}\sum\limits_{-N\leq n\leq N}\frac{1}{2n+1-2s},
\end{align*}
 elementary calculations show that
\begin{align*}
&\sum\limits_{m,n\in\Z,m\neq 0}\frac{f(m)}{(m+ani)^k}(-1)^n=\frac{(-1)^{k-1}i^k}{a^k(k-1)!}\sum\limits_{m=1}^\infty\left(f(m)+(-1)^kf(-m)\right)\frac{d^{k-1}}{ds^{k-1}}\left(\frac{\pi}{\sin(\pi s)}\right)_{s=mi/a},\\&
\sum\limits_{m,n\in\Z,m\neq 0}\frac{g(m)}{(m+ani)^p}=\frac{(-1)^{p-1}i^p}{a^p(p-1)!}\sum\limits_{m=1}^\infty\left(g(m)+(-1)^pg(-m)\right)\frac{d^{p-1}}{ds^{p-1}}\left(\pi \cot(\pi s)\right)_{s=mi/a},\\
&\sum\limits_{m,n\in\Z,m\neq 0} \frac{f(m)(-1)^n}{(m+a(2n+1)i)^k}=\frac{(-1)^ki^k}{(2a)^k(k-1)!}\sum\limits_{m=1}^\infty \left(f(m)+(-1)^{k-1}f(-m)\right)\frac{d^{k-1}}{ds^{k-1}}\left(\frac{\pi}{\cos(\pi s)}\right)_{s=mi/2a},\\&
\sum\limits_{m,n\in\Z,m\neq 0} \frac{g(m)}{(m+a(2n+1)i)^p}=\frac{(-1)^pi^p}{(2a)^p(p-1)!}\sum\limits_{m=1}^\infty \left(g(m)+(-1)^{p}g(-m)\right)\frac{d^{p-1}}{ds^{p-1}}\left(\pi \tan(\pi s)\right)_{s=mi/2a}.
\end{align*}
In general, we have
\begin{align}
&\sum\limits_{m,n\in\Z,m\neq 0}\frac{f(m)(-1)^n}{(bm+c+ani)^k}=\frac{(-1)^{k-1}i^k}{a^k(k-1)!}\sum\limits_{m=1}^\infty \left\{\begin{array}{l} f(m)\frac{d^{k-1}}{ds^{k-1}}\left(\frac{\pi}{\sin(\pi s)}\right)_{s=(bm+c)i/a}\\ +(-1)^kf(-m)\frac{d^{k-1}}{ds^{k-1}}\left(\frac{\pi}{\sin(\pi s)}\right)_{s=(bm-c)i/a}\end{array} \right\},\label{3.1}\\
&\sum\limits_{m,n\in\Z,m\neq 0}\frac{g(m)}{(bm+c+ani)^p}=\frac{(-1)^{p-1}i^p}{a^p(p-1)!}\sum\limits_{m=1}^\infty \left\{\begin{array}{l} g(m)\frac{d^{p-1}}{ds^{p-1}}\left(\pi \cot(\pi s)\right)_{s=(bm+c)i/a}\\ +(-1)^p g(-m)\frac{d^{p-1}}{ds^{p-1}}\left(\pi \cot(\pi s)\right)_{s=(bm-c)i/a}\end{array} \right\},\label{3.12}\\
&\sum\limits_{m,n\in\Z,m\neq 0}\frac{f(m)(-1)^n}{(bm+c+a(2n+1)i)^k}=\frac{(-1)^ki^k}{(2a)^k(k-1)!}\nonumber\\&\quad\quad\quad\quad\quad\quad\quad\quad\quad\quad\quad\quad\quad\quad\times\sum\limits_{m=1}^\infty \left\{\begin{array}{l} f(m)\frac{d^{k1}}{ds^{k-1}}\left(\frac{\pi}{\cos(\pi s)}\right)_{s=(bm+c)i/2a}\\ +(-1)^{k-1}f(-m)\frac{d^{k-1}}{ds^{k-1}}\left(\frac{\pi}{\cos(\pi s)}\right)_{s=(bm-c)i/2a} \end{array} \right\},\label{3.3}\\
&\sum\limits_{m,n\in\Z,m\neq 0}\frac{g(m)}{(bm+c+a(2n+1)i)^p}=\frac{(-1)^pi^p}{(2a)^p(p-1)!}\nonumber\\&\quad\quad\quad\quad\quad\quad\quad\quad\quad\quad\quad\quad\quad\quad\times\sum\limits_{m=1}^\infty \left\{\begin{array}{l} g(m)\frac{d^{p-1}}{ds^{p-1}}(\pi \tan(\pi s))_{s=(bm+c)i/2a}\\ +(-1)^pg(-m)\frac{d^{p-1}}{ds^{p-1}}(\pi \tan(\pi s))_{s=(bm-c)i/2a} \end{array} \right\},\label{3.4}
\end{align}
where $a,b\in \mathbb{R}\setminus \{0\}$ and $c\in \mathbb{R}$.
Then with the help of Theorems \ref{2.2}, \ref{thm2.3} and Corollaries \ref{cor2.4}, \ref{cor2.5}, we can get the following theorems.

\begin{thm}\label{thm3.1} For positive integer $k$ and real $a\in \mathbb{R}\setminus \{0\}$, we have
\begin{align}
\sum\limits_{m,n\in\Z,m\neq 0}\frac{f(m)(-1)^n}{(m+ani)^{2k}}=\frac{\pi^{2k}}{a^{2k}(2k-1)!}\sum\limits_{l=0}^{k-1}&(2l+1)!\left\{\sum\limits_{\mid{\bf r}\mid_l=k-1-l}\prod\limits_{h=0}^l (2h+1)^{2r_h}\right\}\nonumber\\&\times\sum\limits_{m=1}^\infty\frac{\left(f(m)+f(-m)\right)\cosh(m\pi/a)}{\sinh^{2l+2}(m\pi/a)},\label{3.5}\\
\sum\limits_{m,n\in\Z,m\neq 0}\frac{f(m)(-1)^n}{(m+ani)^{2k-1}}=\frac{\pi^{2k-1}}{a^{2k-1}(2k-2)!}\sum\limits_{l=0}^{k-1}&(2l)!\left\{\sum\limits_{\mid{\bf r}\mid_l=k-1-l}\prod\limits_{h=0}^l (2h+1)^{2r_h}\right\}\nonumber\\&\times\sum\limits_{m=1}^\infty\frac{\left(f(m)-f(-m)\right)}{\sinh^{2l+1}(m\pi/a)}.
\end{align}
\end{thm}

\begin{thm}\label{thm3.2}  For positive integer $k$ and real $a\in \mathbb{R}\setminus \{0\}$, we have
\begin{align}
\sum\limits_{m,n\in\Z,m\neq 0}\frac{g(m)}{(m+ani)^{2k}}=&\frac{\pi^{2k}}{a^{2k}(2k-1)!}\sum\limits_{1\leq l\leq j \leq k} \left\{ \binom{2k}{2j}(2l-1)!-\binom{2k}{2j+1}(2l-1)!(2l+1) \right\}\nonumber\\&\quad\quad\quad\quad\times \left\{\sum\limits_{\mid{\bf r}\mid_l=j-l}\prod\limits_{h=0}^l (2h+1)^{2r_h}\right\}\sum\limits_{m=1}^\infty\frac{g(m)+g(-m)}{\sinh^{2l}(m \pi/a)},\label{3.7}
\end{align}
\begin{align}
\sum\limits_{m,n\in\Z,m\neq 0}\frac{g(m)}{(m+ani)^{2k+1}}=&\frac{\pi^{2k+1}}{a^{2k+1}(2k)!}\sum\limits_{0\leq l\leq j \leq k} \left\{ \binom{2k}{2j}(2l)!-\binom{2k}{2j+1}(2l+1)! \right\}\nonumber\\&\quad\quad\times \left\{\sum\limits_{\mid{\bf r}\mid_l=j-l}\prod\limits_{h=0}^l (2h+1)^{2r_h}\right\}\sum\limits_{m=1}^\infty\frac{(g(m)-g(-m))\cosh(m\pi/a)}{\sinh^{2l+1}(m \pi/a)}.
\end{align}
\end{thm}

\begin{thm}\label{thm3.3}  For positive integer $k$ and real $a\in \mathbb{R}\setminus \{0\}$, we have
\begin{align}
\sum\limits_{m,n\in\Z,m\neq 0}\frac{f(m)(-1)^n}{(m+a(2n+1)i)^{2k}}=\frac{\pi^{2k}i}{(2a)^{2k}(2k-1)!}\sum\limits_{l=0}^{k-1}&(-1)^{l-1}(2l+1)!\left\{\sum\limits_{\mid{\bf r}\mid_l=k-1-l}\prod\limits_{h=0}^l (2h+1)^{2r_h}\right\}\nonumber\\&\times\sum\limits_{m=1}^\infty\frac{\left(f(m)-f(-m)\right)\sinh(m\pi/2a)}{\cosh^{2l+2}(m\pi/2a)},
\end{align}
\begin{align}
\sum\limits_{m,n\in\Z,m\neq 0}\frac{f(m)(-1)^n}{(m+a(2n+1)i)^{2k-1}}=\frac{\pi^{2k-1}i}{(2a)^{2k-1}(2k-2)!}\sum\limits_{l=0}^{k-1}&(-1)^{l-1}(2l)!\left\{\sum\limits_{\mid{\bf r}\mid_l=k-1-l}\prod\limits_{h=0}^l (2h+1)^{2r_h}\right\}\nonumber\\&\times\sum\limits_{m=1}^\infty\frac{\left(f(m)+f(-m)\right)}{\cosh^{2l+1}(m\pi/2a)}.
\end{align}
\end{thm}

\begin{thm}\label{thm3.4}  For positive integer $k$ and real $a\in \mathbb{R}\setminus \{0\}$, we have
\begin{align}
\sum\limits_{m,n\in\Z,m\neq 0}\frac{g(m)}{(m+a(2n+1)i)^{2k}}=&\frac{\pi^{2k}}{(2a)^{2k}(2k-1)!}\sum\limits_{1\leq l\leq j \leq k} \left\{ \binom{2k}{2j}(2l-1)!-\binom{2k}{2j+1}(2l-1)!(2l+1) \right\}\nonumber\\&\quad\quad\quad\quad\times (-1)^l \left\{\sum\limits_{\mid{\bf r}\mid_l=j-l}\prod\limits_{h=0}^l (2h+1)^{2r_h}\right\}\sum\limits_{m=1}^\infty\frac{g(m)+g(-m)}{\cosh^{2l}(m \pi/a)},
\end{align}
\begin{align}
\sum\limits_{m,n\in\Z,m\neq 0}\frac{g(m)}{(m+a(2n+1)i)^{2k+1}}=&\frac{\pi^{2k+1}}{(2a)^{2k+1}(2k)!}\sum\limits_{0\leq l\leq j \leq k} \left\{ \binom{2k}{2j}(2l)!-\binom{2k}{2j+1}(2l+1)! \right\}(-1)^l\nonumber\\&\times \left\{\sum\limits_{\mid{\bf r}\mid_l=j-l}\prod\limits_{h=0}^l (2h+1)^{2r_h}\right\}\sum\limits_{m=1}^\infty\frac{(g(m)-g(-m))\sinh(m\pi/a)}{\cosh^{2l+1}(m \pi/a)}.
\end{align}
\end{thm}

From Theorems \ref{thm3.1}-\ref{thm3.4}, we establish many relations of level 2 Eisenstein series. For instance, let
\begin{align}
a_{k,l}:=(2l-1)\frac{\pi^{2k}(-1)^{k+l}}{a^{2k}(2k-1)!}A_{k-1,1}(l-1),\quad k,l\in\N
\end{align}
and
\begin{align}
b_{k,l}:=\frac{\pi^{2k}(-1)^{l}}{a^{2k}(2k-1)!}\sum\limits_{j=l}^k \left\{\frac{1}{2l}\binom{2k}{2j}-\binom{2k}{2j+1}\frac{2l+1}{2l}\right\}(-1)^jA_{j,1}(l)  ,\quad k,l\in\N.
\end{align}
Then, the formulas (\ref{3.5}) and (\ref{3.7}) can be rewritten as
\begin{align}
&\sum\limits_{m,n\in\Z,m\neq 0}\frac{f(m)(-1)^n}{(m+ani)^{2k}}=\sum\limits_{l=1}^{k}a_{k,l}\sum\limits_{m=1}^\infty\frac{\left(f(m)+f(-m)\right)\cosh(m\pi/a)}{\sinh^{2l}(m\pi/a)},\\
&\sum\limits_{m,n\in\Z,m\neq 0}\frac{g(m)}{(m+ani)^{2k}}=\sum\limits_{l=1}^{k}b_{k,l}\sum\limits_{m=1}^\infty\frac{g(m)+g(-m)}{\sinh^{2l}(m\pi/a)}.
\end{align}
Define two square Matrix ${\bf A}_k$ and ${\bf B}_k$ by
\begin{align}
{\bf A}_k:=\{a_{i,j}\}_{k\times k}\quad {\rm and} \quad {\bf B}_k:=\{b_{i,j}\}_{k\times k},
\end{align}
and let
\begin{align}
&{F}_k(f(\cdot),a):=\sum\limits_{m,n\in\Z,m\neq 0}\frac{f(m)(-1)^n}{(m+ani)^{2k}}, \quad  G_k(g(\cdot),a):=\sum\limits_{m,n\in\Z,m\neq 0}\frac{g(m)}{(m+ani)^{2k}},\\
&{\bar F}_l(f(\cdot),a):=\sum\limits_{m=1}^\infty\frac{\left(f(m)+f(-m)\right)\cosh(m\pi/a)}{\sinh^{2l}(m\pi/a)},\quad {\bar G}_l(g(\cdot),a):=\sum\limits_{m=1}^\infty\frac{g(m)+g(-m)}{\sinh^{2l}(m\pi/a)}.
\end{align}
then we have
\begin{align}\label{3.20}
{\bf F}_k(f(\cdot),a)={\bf A}_k\cdot {\bf {\bar F}}_k(f(\cdot),a) \quad {\rm and}\quad {\bf G}_k(g(\cdot),a)={\bf B}_k\cdot {\bf {\bar G}_k}(g(\cdot),a),
\end{align}
where
\begin{align*}
&{\bf F}_k(f(\cdot),a):=({F}_1(f(\cdot),a),{F}_2(f(\cdot),a),\cdots,{F}_k(f(\cdot),a))^{\bf T},\\
&{\bf {\bar F}}_k (f(\cdot),a):=({\bar F}_1(f(\cdot),a),{\bar F}_2(f(\cdot),a),\cdots,{\bar F}_k(f(\cdot),a))^{\bf T},\\
&{\bf G}_k(g(\cdot),a):=({G}_1(g(\cdot),a),{G}_2(g(\cdot),a),\cdots,{G}_k(g(\cdot),a))^{\bf T},\\
&{\bf {\bar G}}_k(g(\cdot),a) :=({\bar G}_1(g(\cdot),a),{\bar G}_2(g(\cdot),a),\cdots,{\bar G}_k(g(\cdot),a))^{\bf T},
\end{align*}
where ${\bf A}^{\bf T}$ is transposed matrix of ${\bf A}$.
We note that if $f(m)=g(m)/\cosh(m\pi/a)$ then
\begin{align*}
{\bar F}_l\left(\frac{g(\cdot)}{\cosh((\cdot)\pi/a)},a\right)={\bar G}_l(g(\cdot),a)\quad {\rm and}\quad {\bf {\bar F}}_l\left(\frac{g(\cdot)}{\cosh((\cdot)\pi/a)},a\right)={\bf {\bar G}}_l(g(\cdot),a)
\end{align*}
Hence, by (\ref{3.20}) one obtain
\begin{align}\label{3.21}
{\bf A}_k^{-1}{\bf F}_k\left(\frac{g(\cdot)}{\cosh((\cdot)\pi/a)},a\right)={\bf B}_k^{-1}{\bf G}_k(g(\cdot),a),
\end{align}
where ${\bf A}^{-1}$ is inversion of matrix ${\bf A}$.
From (\ref{3.21}) we obtain the relations between ${F}_k\left(\frac{g(\cdot)}{\cosh((\cdot)\pi/a)},a\right)$ and $G_k(g(\cdot),a)$. For example,
\begin{align*}
&\sum\limits_{m,n\in\Z,m\neq 0}\frac{g(m)(-1)^n}{\cosh(m\pi/a)(m+ani)^{2}}=\sum\limits_{m,n\in\Z,m\neq 0}\frac{g(m)}{(m+ani)^{2}},\\
&\sum\limits_{m,n\in\Z,m\neq 0}\frac{g(m)(-1)^n}{\cosh(m\pi/a)(m+ani)^{4}}=\sum\limits_{m,n\in\Z,m\neq 0}\frac{g(m)}{(m+ani)^{4}}-\frac{\pi^2}{2a^2}\sum\limits_{m,n\in\Z,m\neq 0}\frac{g(m)}{(m+ani)^{2}},\\
&\sum\limits_{m,n\in\Z,m\neq 0}\frac{g(m)(-1)^n}{\cosh(m\pi/a)(m+ani)^{4}}+\frac{\pi^2}{2a^2}\sum\limits_{m,n\in\Z,m\neq 0}\frac{g(m)(-1)^n}{\cosh(m\pi/a)(m+ani)^{2}}=\sum\limits_{m,n\in\Z,m\neq 0}\frac{g(m)}{(m+ani)^{4}}.
\end{align*}

\subsection{Corollaries}

From Theorems \ref{thm3.1}-\ref{thm3.4} we give the following corollaries.

\begin{cor}\label{cor3.5} For real $a\in \mathbb{R}\setminus \{0\}$, we have
\begin{align}
&\sum\limits_{m,n\in\Z,m\neq 0}\frac{f(m)(-1)^n}{m+ani}= \frac{\pi}{a}\sum\limits_{m=1}^\infty \frac{f(m)-f(-m)}{\sinh(m\pi/a)},\\
&\sum\limits_{m,n\in\Z,m\neq 0}\frac{f(m)(-1)^n}{(m+ani)^{2}}= \frac{\pi^2}{a^2}\sum\limits_{m=1}^\infty \frac{(f(m)+f(-m))\cosh(m\pi/a)}{\sinh^2(m\pi/a)},\\
&\sum\limits_{m,n\in\Z,m\neq 0}\frac{f(m)(-1)^n}{(m+ani)^{3}}= \frac{\pi^3}{2a^3}\sum\limits_{m=1}^\infty \left(\frac{f(m)-f(-m)}{\sinh(m\pi/a)}+2\frac{f(m)-f(-m)}{\sinh^3(m\pi/a)}\right),\\
&\sum\limits_{m,n\in\Z,m\neq 0}\frac{f(m)(-1)^n}{(m+ani)^{4}}= \frac{\pi^4}{6a^4}\sum\limits_{m=1}^\infty (f(m)+f(-m))\left(\frac{\cosh(m\pi/a)}{\sinh^2(m\pi/a)}+6\frac{\cosh(m\pi/a)}{\sinh^4(m\pi/a)}\right).
\end{align}
\end{cor}

\begin{cor} \label{cor3.6} For real $a\in \mathbb{R}\setminus \{0\}$, we have
\begin{align}
&\sum\limits_{m,n\in\Z,m\neq 0}\frac{g(m)}{(m+ani)^{2}}=\frac{\pi^2}{a^2}\sum\limits_{m=1}^\infty \frac{g(m)+g(-m)}{\sinh^2(m\pi/a)},\\
&\sum\limits_{m,n\in\Z,m\neq 0}\frac{g(m)}{(m+ani)^{3}}=\frac{\pi^3}{a^3}\sum\limits_{m=1}^\infty \frac{(g(m)-g(-m))\cosh(m\pi/a)}{\sinh^3(m\pi/a)},\\
&\sum\limits_{m,n\in\Z,m\neq 0}\frac{g(m)}{(m+ani)^{4}}=\frac{2\pi^4}{3a^4}\sum\limits_{m=1}^\infty \frac{g(m)+g(-m)}{\sinh^2(m\pi/a)}+\frac{\pi^4}{a^4}\sum\limits_{m=1}^\infty \frac{g(m)+g(-m)}{\sinh^4(m\pi/a)}.
\end{align}
\end{cor}

\begin{cor}\label{cor3.7} For real $a\in \mathbb{R}\setminus \{0\}$, we have
\begin{align}
&\sum\limits_{m,n\in\Z,m\neq 0}\frac{f(m)(-1)^n}{m+a(2n+1)i}=-\frac{\pi i}{2a}\sum\limits_{m=1}^\infty \frac{f(m)+f(-m)}{\cosh(m\pi/2a)},\\
&\sum\limits_{m,n\in\Z,m\neq 0}\frac{f(m)(-1)^n}{(m+a(2n+1)i)^2}=-\frac{\pi^2 i}{4a^2}\sum\limits_{m=1}^\infty \frac{(f(m)-f(-m))\sinh(m\pi/2a)}{\cosh^2(m\pi/2a)},\\
&\sum\limits_{m,n\in\Z,m\neq 0}\frac{f(m)(-1)^n}{(m+a(2n+1)i)^3}=-\frac{\pi^3 i}{16a^2}\sum\limits_{m=1}^\infty \left(\frac{f(m)+f(-m)}{\cosh(m\pi/2a)}-2\frac{f(m)+f(-m)}{\cosh^3(m\pi/2a)}\right),\\
&\sum\limits_{m,n\in\Z,m\neq 0}\frac{f(m)(-1)^n}{(m+a(2n+1)i)^4}=-\frac{\pi^4 i}{96a^4}\sum\limits_{m=1}^\infty (f(m)-f(-m))\left(\frac{\sinh(m\pi/2a)}{\cosh^2(m\pi/2a)}-6\frac{\sinh(m\pi/2a)}{\cosh^4(m\pi/2a)}\right).
\end{align}
\end{cor}

\begin{cor}\label{cor3.8} For real $a\in \mathbb{R}\setminus \{0\}$, we have
\begin{align}
&\sum\limits_{m,n\in\Z,m\neq 0}\frac{g(m)}{(m+a(2n+1)i)^2}=-\frac{\pi^2}{4a^2}\sum\limits_{m=1}^\infty \frac{g(m)+g(-m)}{\cosh^2(m\pi/2a)},\\
&\sum\limits_{m,n\in\Z,m\neq 0}\frac{g(m)}{(m+a(2n+1)i)^3}=-\frac{\pi^3}{8a^3}\sum\limits_{m=1}^\infty \frac{(g(m)-g(-m))\sinh(m\pi/2a)}{\cosh^3(m\pi/2a)},\\
&\sum\limits_{m,n\in\Z,m\neq 0}\frac{g(m)}{(m+a(2n+1)i)^4}=-\frac{\pi^4}{48a^4}\sum\limits_{m=1}^\infty \left(2\frac{g(m)+g(-m)}{\cosh^2(m\pi/2a)} -3\frac{g(m)+g(-m)}{\cosh^4(m\pi/2a)}\right).
\end{align}
\end{cor}

\subsection{Examples}

Since the four infinite series involving hyperbolic functions
\begin{align*}
\sum\limits_{n=1}^\infty \frac{1}{\sinh^{2m}(n\pi)},\ \sum\limits_{n=1}^\infty \frac{1}{\cosh^{2m}(n\pi)},\ \sum\limits_{n=1}^\infty \frac{n^{2p}}{\sinh^{2k}(n\pi)}\quad {\rm and}\quad \sum\limits_{n=1}^\infty \frac{n^{2p}}{\cosh^{2k}(n\pi)}\quad (m,p,k\in\N, p\geq k)
\end{align*}
can be evaluated by Gamma function and $\pi$ (explicit evaluations see \cite{L1974,C2018}). For example, we have
\begin{align*}
&\sum\limits_{n=1}^\infty \frac{1}{\sinh^2(n\pi)}=\frac{1}{6}-\frac{1}{2\pi},\\
&\sum\limits_{n=1}^\infty \frac{1}{\sinh^4(n\pi)}=-\frac{11}{90}+\frac{1}{3\pi}+\frac{\Gamma^8(1/4)}{960\pi^2},\\
&\sum\limits_{n=1}^\infty \frac{1}{\cosh^2(n\pi)}=-\frac 1{2}+\frac 1{2\pi}+\frac{\Gamma^4(1/4)}{16\pi^3},\\
&\sum\limits_{n=1}^\infty \frac{1}{\cosh^4(n\pi)}=-\frac{1}{2}+\frac{1}{3\pi}+\frac{\Gamma^4(1/4)}{24\pi^3}+\frac{\Gamma^8(1/4)}{192\pi^6}.
\end{align*}
Hence, from Corollaries \ref{cor3.5}-\ref{cor3.8}, we give some well-known and new results of level 2 Eisenstein series involving hyperbolic functions
\begin{align*}
&\sum\limits_{m,n\in\Z,m\neq 0}\frac{m^2(-1)^n}{\sinh(m\pi)(m+ni)}=-\frac 1{4\pi}+\frac{\Gamma^8(1/4)}{768\pi^5},\\
&\sum\limits_{m,n\in\Z,m\neq 0}\frac{m^4(-1)^n}{\sinh(m\pi)(m+ni)}=\frac{\Gamma^8(1/4)}{640\pi^6},\\
&\sum\limits_{m,n\in\Z,m\neq 0}\frac{m^6(-1)^n}{\sinh(m\pi)(m+ni)}=\frac{\Gamma^{16}(1/4)}{57344\pi^{11}},\\
&\sum\limits_{m,n\in\Z,m\neq 0}\frac{m^8(-1)^n}{\sinh(m\pi)(m+ni)}=\frac{3\Gamma^{16}(1/4)}{40960\pi^{12}},\\
&\sum\limits_{m,n\in\Z,m\neq 0}\frac{m^2(-1)^n}{\cosh(m\pi)(m+(2n+1)i/2)}=i\left(\frac{\Gamma^4(1/4)}{32\pi^3}-\frac{\Gamma^8(1/4)}{768\pi^5}\right),\\
&\sum\limits_{m,n\in\Z,m\neq 0}\frac{m^4(-1)^n}{\cosh(m\pi)(m+(2n+1)i/2)}=i\left(\frac{3\Gamma^8(1/4)}{5120\pi^6}-\frac{\Gamma^{12}(1/4)}{8192\pi^8}\right),\\
&\sum\limits_{m,n\in\Z,m\neq 0}\frac{m^6(-1)^n}{\cosh(m\pi)(m+(2n+1)i/2)}=-i\left(\frac{9\Gamma^{12}(1/4)}{32768\pi^9}+\frac{3\Gamma^{16}(1/4)}{1835008\pi^{11}}\right),\\
&\sum\limits_{m,n\in\Z,m\neq 0}\frac{m^8(-1)^n}{\cosh(m\pi)(m+(2n+1)i/2)}=-i\left(\frac{21\Gamma^{16}(1/4)}{5242880\pi^{12}}+\frac{33\Gamma^{20}(1/4)}{8388608\pi^{14}}\right),\\
&\sum\limits_{m,n\in\Z,m\neq 0} \frac{m^2}{\cosh^2(m\pi/2)(m+2ni)^2}=-\frac 1{4}+\frac{\Gamma^8(1/4)}{768\pi^4},\\
&\sum\limits_{m,n\in\Z,m\neq 0} \frac{1}{\sinh^2(m\pi)(m+ni)^2}=-\frac{11}{45}\pi^2+\frac 2{3}\pi+\frac{\Gamma^8(1/4)}{960\pi^4},\\
&\sum\limits_{m,n\in\Z,m\neq 0} \frac{1}{\sinh(2m\pi)(m+ni)^3}=-\frac{11}{90}\pi^2+\frac 1{3}\pi^2+\frac{\Gamma^8(1/4)}{1920\pi^4},\\
&\sum\limits_{m,n\in\Z,m\neq 0}\frac{(-1)^n}{\sinh(m\pi)\left(m+(2n+1)i/2\right)^2}=i\left(\pi^2-\pi-\frac{\Gamma^4(1/4)}{8\pi}\right),\\
&\sum\limits_{m,n\in\Z,m\neq 0}\frac{(-1)^n}{\sinh(m\pi)\left(m+(2n+1)i/2\right)^4}=i\left(\frac {\pi^3}{2}-\frac{5}{6}\pi^4+\frac{\pi}{16}\Gamma^4(1/4)+\frac{\Gamma^8(1/4)}{96\pi^2}\right),\\
&\sum\limits_{m,n\in\Z,m\neq 0}\frac{(-1)^n}{\cosh(m\pi)\left(m+(2n+1)i/2\right)}=i\left(\pi-1-\frac{\Gamma^4(1/4)}{8\pi^2}\right),\\
&\sum\limits_{m,n\in\Z,m\neq 0}\frac{(-1)^n}{\cosh(m\pi)\left(m+(2n+1)i/2\right)^3}=i\left(\frac {\pi^2}{6}-\frac{\pi^3}{2}+\frac{\Gamma^4(1/4)}{48}+\frac{\Gamma^8(1/4)}{96\pi^3}\right),\\
&\sum\limits_{m,n\in\Z,m\neq 0}\frac{1}{\sinh^2(m\pi/2)\big(m+(2n+1)i\big)^2}=\pi-\frac{\pi^2}{3},\\
&\sum\limits_{m,n\in\Z,m\neq 0}\frac{1}{\cosh^2(m\pi)\left(m+(2n+1)i/2\right)^4}=\frac{4}{45}\pi^3-\frac{\pi^4}{3}+\frac{\pi}{90}\Gamma^4(1/4)+\frac{\Gamma^8(1/4)}{288\pi^2}+\frac{\Gamma^{12}(1/4)}{1280\pi^5},\\
&\sum\limits_{m,n\in\Z,m\neq 0}\frac{1}{\sinh(2m\pi)\big(m+(2n+1)i/2\big)^3}=\frac{\pi^3}{2}-\frac{\pi^2}{3}-\frac{\Gamma^4(1/4)}{24}-\frac{\Gamma^8(1/4)}{192\pi^3}.
\end{align*}
The tenth equation  appear as example of Example 3 in the \cite{T2015}. It should be emphasized that the reference \cite{T2015} also contains many other types of double Eisenstein series.

It is possible that many other evaluations of double Eisenstein series involving hyperbolic functions
can be obtained by using the methods and techniques of the present paper. For example, by Theorems \ref{thm3.1}-\ref{thm3.4}, it is clear that the twelve level 2 Eisenstein series involving hyperbolic functions
\begin{align*}
&\sum\limits_{m,n\in\Z,m\neq 0} \frac{\coth^{2p}(m\pi)}{(m+ni)^{2k+2}},\\
&\sum\limits_{m,n\in\Z,m\neq 0}  \frac{\coth^{2p+1}(m\pi)}{(m+ni)^{2k+1}},\\
&\sum\limits_{m,n\in\Z,m\neq 0} \frac{\tanh^{2p}(m\pi)}{(m+(2n+1)i/2)^{2k+2}},\\
&\sum\limits_{m,n\in\Z,m\neq 0} \frac{\tanh^{2p+1}(m\pi)}{(m+(2n+1)i/2)^{2k+1}},\\
&\sum\limits_{m,n\in\Z,m\neq 0}\frac{(-1)^n}{(m+ni)^{2k}\cosh(m\pi)\sinh^{2p}(m\pi)},\\
& \sum\limits_{m,n\in\Z,m\neq 0}\frac{(-1)^n}{(m+ni)^{2k-1}\sinh^{2p+1}(m\pi)},\\
&\sum\limits_{m,n\in\Z,m\neq 0}\frac{1}{(m+ni)^{2k+1}\cosh(m\pi)\sinh^{2p+1}(m\pi)},\\
& \sum\limits_{m,n\in\Z,m\neq 0}\frac{1}{(m+ni)^{2k}\sinh^{2p+2}(m\pi)},\\
&\sum\limits_{m,n\in\Z,m\neq 0}\frac{(-1)^n}{(m+(2n+1)i/2)^{2k}\sinh(m\pi)\cosh^{2p}(m\pi)},\\
& \sum\limits_{m,n\in\Z,m\neq 0}\frac{(-1)^n}{(m+(2n+1)i/2)^{2k-1}\cosh^{2p+1}(m\pi)},\\
&\sum\limits_{m,n\in\Z,m\neq 0}\frac{1}{(m+(2n+1)i/2)^{2k+1}\sinh(m\pi)\cosh^{2p+1}(m\pi)},\\
&\sum\limits_{m,n\in\Z,m\neq 0}\frac{1}{(m+(2n+1)i/2)^{2k}\cosh^{2p+2}(m\pi)}.
\end{align*}
can be represented by $\sum\limits_{n=1}^\infty \frac{1}{\sinh^{2l}(n\pi)}$ or $\sum\limits_{n=1}^\infty \frac{1}{\cosh^{2l}(n\pi)}$, which implies that these level 2 Eisenstein series
can be expressed in terms of $\Gamma$ function and $\pi$. Here $k\in\N$ and $p\in\N_0$. For example, the evaluations of the first and second sums can be obtained by Theorem \ref{thm3.2}. The evaluations of the third and forth sums can be deduced by Theorem \ref{thm3.4}.

Similarly, by Theorems \ref{thm2.2}, \ref{thm2.3}, Corollaries \ref{cor2.4}, \ref{cor2.5},  and formulas (\ref{3.1})-(\ref{3.4}) we can give explicit evaluations of the sums 
\begin{align*}
&\sum\limits_{m,n\in\Z} \frac{f(m)(-1)^n}{(2bm+b+ani)^k},\ \sum\limits_{m,n\in\Z} \frac{g(m)}{(2bm+b+ani)^p},\\
&\sum\limits_{m,n\in\Z} \frac{f(m)(-1)^n}{(2bm+b+a(2n+1)i)^k},\ \sum\limits_{m,n\in\Z} \frac{g(m)}{(2bm+b+a(2n+1)i)^p},
\end{align*}
where $a,b\in \mathbb{R}\setminus \{0\}$. For instance, a simple calculation gives
\begin{align}
\sum\limits_{m,n\in\Z} \frac{f(m)(-1)^n}{(2bm+b+ani)^{2k}}=&\frac{\pi^{2k}}{a^{2k}(2k-1)!}\sum\limits_{l=0}^{k-1}(2l+1)!\left\{\sum\limits_{\mid{\bf r}\mid_l=k-1-l}\prod\limits_{h=0}^l (2h+1)^{2r_h}\right\}\nonumber\\ &\quad\quad\quad\quad\quad\times \sum\limits_{m=1}^\infty \frac{(f(m-1)+f(-m))\cosh((2m-1)b\pi/a)}{\sinh^{2l+2}((2m-1)b\pi/a)}.\label{3.26}
\end{align}
Setting $k=1$ in equation above yields
\begin{align}
\sum\limits_{m,n\in\Z} \frac{f(m)(-1)^n}{(2bm+b+ani)^{2}}=\frac{\pi^2}{a^2} \sum\limits_{m=1}^\infty \frac{(f(m-1)+f(-m))\cosh((2m-1)b\pi/a)}{\sinh^{2}((2m-1)b\pi/a)}.
\end{align}
Hence, we can deduce the two evaluations
\begin{align*}
&\sum\limits_{m,n\in\Z} \frac{(-1)^n}{\cosh((2m+1)\pi/2)(2m+1+2ni)^{2}}=-\frac{\pi}{4}+\frac{\Gamma^4(1/4)}{32\pi},\\
&\sum\limits_{m,n\in\Z} \frac{(-1)^n}{\cosh^3((2m+1)\pi/2)(2m+1+2ni)^{2}}=-\frac{\pi}{2}+\frac{\Gamma^4(1/4)}{32\pi},
\end{align*}
where we used the three well known results (see \cite{L1974,L1975})
\begin{align*}
& \su \frac{1}{\sinh^2(n\pi)}=\frac 1{6}-\frac{1}{2\pi},\\
&\su \frac{(-1)^{n-1}}{\sinh^2(n\pi)}=-\frac 1{6}+\frac{\Gamma^4(1/4)}{32\pi^3},\\
&\su \frac{1}{\sinh^2((2n-1)\pi/2)}=-\frac{1}{2\pi}+\frac{\Gamma^4(1/4)}{16\pi^3}.
\end{align*}
Note that in \cite{L1974,L1975},
$u=\frac{\Gamma^4(1/4)}{8\pi}$. It is obvious that the results of Tsumura \cite{T2008,T2009,T2010} can be established by using the methods of the present paper.

We finally remark that using the method of this paper, it is possible to evaluate other hyperbolic series. For example, in recent paper \cite{T2018}, Tsumura studied the double series 
$\tilde{H}_{2k}(\tau)$ ($\tau$ belongs to the upper half-plane) and proved that $\tilde{H}_{2k}(i)\in \mathbb{Q}[\pi, \Gamma^2(1/4)/(2\sqrt{2\pi})]$, where $\tilde{H}_{2k}(\tau)$ is defined by the double series (Though we treats only the case $\tau=i$ in this paper, it is desirable and possible to treat the general case of $\tau$)
\begin{align*}
\tilde{H}_{2k}(\tau)=\sum\limits_{\substack{ m,n\in\Z,\\(m,n)\neq (0,0)}}\frac{(-1)^m}{\cosh(m\pi i/\tau)(m+n\tau)^{2k}}\quad (k\in\N).
\end{align*}
If letting $g(m)=(-1)^m/\cosh(m\pi)$ and $a=1$ in (\ref{3.7}), we find that $\tilde{H}_{2k}(i)$ can be expressed in terms of hyperbolic series 
\begin{align*}
\sum\limits_{m=1}^\infty \frac{(-1)^m}{\cosh(m\pi)\sinh^{2l}(m\pi)} \quad (l\leq k).
\end{align*}
For instance, 
\begin{align*}
\sum\limits_{m,n\in\Z,m\neq 0}\frac{(-1)^m}{\cosh(m\pi)(m+ni)^{4}}=&\frac{4\pi^4}{3}\sum\limits_{m=1}^\infty \frac{(-1)^m}{\cosh(m\pi)\sinh^2(m\pi)}+2\pi^4\sum\limits_{m=1}^\infty \frac{(-1)^m}{\cosh(m\pi)\sinh^4(m\pi)}\\
=&-\frac{83}{360}\pi^4+\frac{\Gamma^2(1/4)}{6\sqrt{2}}\pi^{5/2}-\frac{\Gamma^8(1/4)}{640\pi^2}.
\end{align*}
Conversely, we obtain the conclusion 
\begin{align*}
\sum\limits_{m=1}^\infty \frac{(-1)^m}{\cosh(m\pi)\sinh^{2l}(m\pi)}\in \mathbb{Q}[\pi, \Gamma^2(1/4)/(2\sqrt{2\pi})].
\end{align*}
With the help of results in \cite{T2018}, we have
\begin{align*}
&\sum\limits_{m=1}^\infty \frac{(-1)^m}{\cosh(m\pi)\sinh^2(m\pi)}=\frac{5}{12}-\frac{\Gamma^2(1/4)}{4\sqrt{2}\pi^{3/2}},\\
&\sum\limits_{m=1}^\infty \frac{(-1)^m}{\cosh(m\pi)\sinh^4(m\pi)}=-\frac{283}{720}+\frac{\Gamma^2(1/4)}{4\sqrt{2}\pi^{3/2}}-\frac{\Gamma^8(1/4)}{1280\pi^6}.
\end{align*}


{\small
\begin{thebibliography}{99}

\bibitem{A1974}
R. Ayoub.  {\sl Euler and the zeta function}. Amer. Math. Monthly, 1974, {\bf 81}: 1067-1086.


\bibitem{B1977}
B.C. Berndt. {\sl Modular transformations ang generalizations of several formulae of Ramanujan}. Rocky Mt. J. Math., 1977, {\bf 7}: 147-189.

\bibitem{B1978}
B.C. Berndt. {\sl Analytic Eisenstein series, theta-functions, and series relations in the spirit of Ramanujan}.
J. Reine Angew. Math., 1978, {\bf 304}: 332--365.


\bibitem{B1989} B.C. Berndt. {\sl Ramanujan's notebooks, part II}. Springer-Verlag, New York, 1989.

\bibitem{B1991} B.C. Berndt. {\sl Ramanujan's notebooks, part III}. Springer-Verlag, New York, 1990.

\bibitem{H1927}
G.H. Hardy, P.V.Seshu Aiyar and B.M. Wilson. {\sl Collected papers of Srinivasa Ramanujan}. Cambridge University Press, 1927.

\bibitem{H1962}
A. Hurwitz. {\sl Mathematische Werke. Bd. I: Funktionentheorie }(German), Herausgegeben von der Abteilung
f$\rm \ddot{u}$r Mathematik und Physik der Eidgen$\rm \ddot{o}$ssischen Technischen Hochschule in Z$\rm \ddot{u}$rich Birkh$\rm \ddot{a}$user, Basel, 1962.

\bibitem{T2015}
Y. Komori, K. Matsumoto and H. Tsumura. {\sl Infinite series involving hyperbolic functions}. Lith. Math. J., 2015, {\bf 55}(1): 102-118.

\bibitem{L1974}
C.B. Ling. {\sl On summation of series of hyperbolic functions}. Siam J. Math. Anal., 1974, {\bf 5}(4): 551-562.

\bibitem{L1975}
C.B. Ling. {\sl On summation of series of hyperbolic functions, II}. Siam J. Math. Anal., 1975, {\bf 6}(1): 129-139.

\bibitem{T2008}
H. Tsumura. {\sl On certain analogues of Eisenstein series and their evaluation formulas of Hurwitz type}. Bull. London Math. Soc., 2008, {\bf 40}: 289-297.

\bibitem{T2009}
H. Tsumura. {\sl Evaluation of certain classes of Eisenstein-type series}. Bull. Aust. Math. Soc., 2009, {\bf 79}: 239-247.

\bibitem{T2010}
H. Tsumura. {\sl Analogues of the Hurwitz formulas for level 2 Eisenstein series}. Results. Math., 2010, {\bf 58}: 365-378.

\bibitem{T2012}
H. Tsumura. {\sl Analogues of level-N Eisenstein series}. Pacific J. Math., 2012, {\bf 255}(2): 489-510.

\bibitem{T2018}
H. Tsumura. {\sl Double series identities arising from Jacobi's identity of the theta function}.  Results. Math., 2018, {\bf 73}(10): 1-10.

\bibitem{C2018}
C. Xu. {\sl Evaluations of infinite series involving reciprocal quadratic hyperbolic functions}. arXiv:1801.07565, 2018.

\end{thebibliography}}
 \end{document}